\documentclass[12pt, reqno, letterpaper]{amsart}

\usepackage{amssymb}  
\usepackage{calc}          
\usepackage{multicol}     

\input cyracc.def
\newfam\cyrfam
\font\tencyr=wncyr10 
\font\tenitcyr=wncyi10 
\def\cyr{\fam\cyrfam\tencyr\cyracc}
\def\itcyr{\fam\cyrfam\tenitcyr\cyracc}

\theoremstyle{plain}
\newtheorem{theorem}{Theorem}[section]
\newtheorem{lemma}[theorem]{Lemma}
\newtheorem{definition}[theorem]{Definition}
\newtheorem{proposition}[theorem]{Proposition}
\newtheorem{corollary}[theorem]{Corollary}
\newtheorem{example}[theorem]{Example}
\theoremstyle{remark}

\newtheorem*{acknowledgement}{Acknowledgement}

\newcommand{\seclabel}[1]{\label{sec:#1}} 
\newcommand{\thmlabel}[1]{\label{thm:#1}} 
\newcommand{\lemlabel}[1]{\label{lem:#1}} 
\newcommand{\corolabel}[1]{\label{coro:#1}} 
\newcommand{\proplabel}[1]{\label{prop:#1}} 
\newcommand{\deflabel}[1]{\label{def:#1}} 
\newcommand{\exlabel}[1]{\label{ex:#1}} 

\newcommand{\secref}[1]{\ref{sec:#1}} 
\newcommand{\thmref}[1]{\ref{thm:#1}} 
\newcommand{\lemref}[1]{\ref{lem:#1}} 
\newcommand{\cororef}[1]{\ref{coro:#1}} 
\newcommand{\propref}[1]{\ref{prop:#1}} 
\newcommand{\defref}[1]{\ref{def:#1}} 
\newcommand{\exref}[1]{\ref{ex:#1}} 
\renewcommand{\eqref}[1]{\ref{eq:#1}} 

\newcommand{\Aut}{\mathrm{Aut}}

\newcommand{\Hom}{\mathrm{Hom}}
\newcommand{\RMlt}{\mathrm{RMlt}}
\newcommand{\LMlt}{\mathrm{LMlt}}

\newcommand{\Syl}{\mathrm{Syl}}
\newcommand{\Hall}{\mathrm{Hall}}

\newcommand\sbl[1]{\langle#1\rangle}   
\newcommand\normal{\trianglelefteq}  
\newcommand\iv{^{-1}} 
\newcommand\ld{\backslash} 
\newcommand\rd{/}  

\newcommand{\RR}{\mathbb{R}}
\newcommand{\ZZ}{\mathbb{Z}}

\newcommand{\TTT}{\mathcal{T}}
\newcommand\res{\mathord {\upharpoonright}}  
\newcommand{\oa}{\overline{\alpha}}

\title[The Structure of Extra Loops]
{The Structure of Extra Loops}

\author[M.~K.~Kinyon]{Michael~K.~Kinyon}
\address{Department of Mathematical Sciences \\
Indiana University South Bend \\
South Bend, IN 46634 USA}
\email{mkinyon@iusb.edu}
\urladdr{http://mypage.iusb.edu/\symbol{126}mkinyon}
\author[K.~Kunen]{Kenneth~Kunen$^*$}
\thanks{$^*$Partially supported by NSF Grant DMS-0097881}
\address{Department of Mathematics \\
University of Wisconsin \\
Madison, WI 57306 USA}
\email{kunen@math.wisc.edu}
\urladdr{http://www.math.wisc.edu/\symbol{126}kunen}

\date{\today}


\begin{document}

\begin{abstract}
The Sylow theorems hold for finite extra loops, as does P. Hall's
theorem for finite solvable extra loops.
Every finite nonassociative extra loop $Q$ has a nontrivial center,
$Z(Q)$.  Furthermore, $Q/Z(Q)$ is a group whenever $|Q| < 512$.
Loop extensions are used to construct an infinite nonassociative
extra loop with a trivial center and a nonassociative
extra loop $Q$ of order $512$ such that $Q/Z(Q)$ is nonassociative.
There are exactly $16$ nonassociative extra loops of order $16p$
for each odd prime $p$.
\end{abstract}

\maketitle

\section{Introduction}
\seclabel{intro}

\begin{definition}
\deflabel{extra}
A loop $Q$ is an \emph{extra loop} iff $Q$
is both conjugacy closed (a CC-loop) and a Moufang loop.
\end{definition}

\begin{lemma}
\lemlabel{eqns}
A loop $Q$ is an extra loop iff $Q$ satisfies
one (equivalently all) of the following equations:
\begin{itemize}
\item[1.]
$(x \cdot y  z) \cdot y = x  y \cdot z  y$.   
\item[2.]
$y  z \cdot y  x = y \cdot (z  y \cdot x)$.   
\item[3.]
$(xy\cdot z) \cdot x = x \cdot (y\cdot zx)$.
\end{itemize}
\end{lemma}

Extra loops were first introduced via these equations 
by Fenyves \cite{FENA, FENB}, who proved the equivalence of (1)(2)(3).
Goodaire and Robinson \cite{GRB} showed that
Definition \defref{extra} is equivalent,
and this definition is often more useful in practice,
since one may combine results in the literature on
CC-loops and on Moufang loops to prove 
theorems about extra loops.

Moufang loops are discussed in standard
texts \cite{Bel, Br, Pf} on loop theory.
In particular, these loops are diassociative by Moufang's Theorem.
CC-loops were introduced by Goodaire and Robinson \cite{GRA, GRB},
and independently (with different terminology) by
{\cyr So\u\i kis} \cite{SO}. Further discussion can be found
in \cite{DRA1, DRA2, KKP, KUNC}.

If $Q$ is an extra loop and $N = N(Q)$ is the nucleus of $Q$,
then $N$ is a normal subloop of $Q$ and
$Q/N$ is a boolean group (see Fenyves \cite{FENB}).
Besides leading to the result of Chein and Robinson that
extra loops are exactly those Moufang loops with squares
in the nucleus \cite{CR}, Fenyves's result  
suggests that one might provide a detailed structure
theory for finite extra loops.  
A start on such a theory was made in \cite{KKP}, where it
was shown that if $Q$ is a finite nonassociative extra loop,
then $|N|$ is even and $|Q : N| \ge 8$, so that 
$16 \mid |Q|$.
The five nonassociative Moufang loops of order 16
are all extra loops
(see Chein \cite{CH}, p.~49).  Among these five is the Cayley loop (1845),
which is the oldest known example of a nonassociative loop.

The Cayley loop is usually described by starting with the 
octonion ring  ($\RR^8$), and restricting the multiplication
to $\{\pm e_i : 0 \le i \le 7\}$, where the $e_i$ are the standard
basis vectors. Restricting to $\RR^8 \backslash \{0\}$ or to $S^7$
does not yield an extra loop (it is Moufang, but not CC).
In fact, by Nagy and Strambach (\cite{NS},
Corollary 2.5, p. 1043),
there are no nonassociative connected smooth extra loops.
There are also no nonassociative connected compact extra loops,
since $Q/N$ is boolean, and hence totally disconnected.

The main results of this paper are listed in the abstract.
After we review basic facts about extra loops in \S\secref{basics},
we characterize the nuclei of nonassociative extra loops in \S\secref{nuc}.
The Sylow theorems are proved in \S\secref{sylow}, and P. Hall's
theorem is proved in \S\secref{Hall}. The center is discussed in
\S\secref{center}. In \S\secref{ext}, we consider loop extensions and
describe the two examples mentioned in the abstract.
In \S\secref{semi} we analyze the nonassociative
extra loops of order $16p$, for $p$ an odd prime,
and show that the number of such loops is independent of $p$;
it follows that this number is $16$, since by \cite{GMR},
there are $16$ such loops of order $48$.


\section{Basic Facts}
\seclabel{basics}

We collect some facts from the literature.  In particular, we
point out that an extra loop yields four boolean groups
which help elucidate the loop structure. One is the quotient
by the nucleus:

\begin{lemma}
\lemlabel{square}
Let $Q$ be an extra loop with nucleus $N = N(Q)$.
\begin{enumerate}
\item[1.] For each $x\in Q$, $x^2 \in N$.
\item[2.] $Q/N$ is a boolean group.
\item[3.] Every finite subloop of $Q$ of odd order is contained in $N$.
\item[4.] Every element of $Q$ of finite odd order is contained in $N$.
\end{enumerate}
\end{lemma}

The lemma, particularly (1), is due to Fenyves \cite{FENB}.
Considered as a Moufang or CC-loop, an extra loop has
a normal nucleus, so (2) follows from (1) and the fact that a Moufang
or CC-loop of exponent 2 is a boolean group. (3) follows
from (2) (since $Q \to Q/N$ maps
the subloop to $\{1\}$), and (4) follows from (3).

\begin{corollary}
\corolabel{Lagrange}
Every finite extra loop has the Lagrange property; that is,
the order of every subloop divides the order of the loop.
\end{corollary}

This follows from the fact that $Q/N$ is a group, so that
both $Q/N$ and $N$ have the Lagrange property;
see Bruck \cite{Br}, \S V.2,  Lemma 2.1.
This corollary holds for all CC-loops $Q$, because
Basarab \cite{BAS} has shown that $Q/N$ is an abelian group;
see also \cite{KKP} for an  exposition of
Basarab's proof, and see \cite{DRA1} for related results.

Another boolean group is generated by the associators:

\begin{definition}
\deflabel{associator}
For each $x,y,z$ in a loop $Q$, define the \emph{associator} $(x,y,z)\in Q$
by $(x\cdot yz) (x,y,z) = xy \cdot z$. Let $A(Q)$ be the subloop of $Q$
generated by all the associators.
\end{definition}

In an extra loop $Q$, $A(Q)\leq N(Q)$, since $Q/N(Q)$ is a group.
Furthermore, by \S5 of \cite{KKP}, we have:

\begin{lemma}
\lemlabel{assoc}
In any extra loop $Q$:
\begin{itemize}
\item[1.] $(x,y,z)$ is invariant under all permutations of the set $\{x,y,z\}$.
\item[2.] $(x,y,z) = (ux,vy,wz)$ for all $x,y,z \in Q$ and $u,v,w \in N(Q)$.
\item[3.] $(x,y,z) = (x\iv, y, z)$.
\item[4.] $(x,y,z)$ commutes with each of $x,y,z$.
\item[5.] $A(Q) \le Z(N(Q))$ and $A(Q)$ is a boolean group.
\end{itemize}
\end{lemma}

Note that Lemma \lemref{assoc} shows that the associator 
$(x,y,z)$ determines a totally symmetric mapping from $(Q/N)^3$ into $A(Q)$.

If $|Q| < 512$, then Theorem \thmref{assoc-512} will show that
$A(Q) \le Z(Q)$ (equivalently, $Q/Z(Q)$ is a group);
this fails for some $Q$ of order $512$; see Example \exref{512}.
For any finite nonassociative extra loop, 
$|Z(Q)\cap A(Q)| \ge 2$ (see Theorem \thmref{center}).

The properties we have listed for associators actually
characterize extra loops:

\begin{lemma}
\lemlabel{extra-assoc}
Suppose that $Q$ is a loop with the following properties:
\begin{itemize}
\item[1.]
 $Q$ is flexible, that is, $(x,y,x) = 1$ for all $x,y\in Q$.
\item[2.]
 Every associator is in the nucleus.
\item[3.]
 The square of every associator is $1$.
\item[4.]
$(x,y,z)$ is invariant under all permutations of \{x,y,z\}.
\item[5.]
$(x,y,z)$ commutes with each of $x,y,z$.
\end{itemize}
Then $Q$ is an extra loop.
\end{lemma}
\begin{proof}
$
x\cdot [y\cdot zx] = x\cdot yz\cdot x\cdot (y,z,x) =
[xy\cdot z]  (x,y,z) x (y,z,x) =
[xy\cdot z] \cdot x .
$
\end{proof}

The third boolean group is the right inner mapping group, which
turns out in this case to coincide with the left inner mapping group
(see \lemref{Mlt}(5) below). We use the following notation.

\begin{definition}
\deflabel{Mlt}
For any loop $Q$, 
the \emph{left translations} $L_x$ and
\emph{right translations} $R_y$ are defined by:
$xy = x R_y = y L_x$. The \emph{right} and
\emph{left multiplication groups} are, respectively
\[
\RMlt = \RMlt(Q) = \sbl{R_y : y\in Q}
\quad \text{and} \quad
\LMlt = \LMlt(Q) = \sbl{L_x : x\in Q} .
\]
For $S\subset Q$, set $R(S) := \{ R_x : x \in S \}$.
The \emph{right} and \emph{left inner mapping groups}
are, respectively,
\begin{align*}
\RMlt_1 &= \RMlt_1(Q) = \{ g \in \RMlt : 1g = 1 \}\ \ \text{and}\\
\LMlt_1 &= \LMlt_1(Q) = \{ g \in \LMlt : 1g = 1\}.
\end{align*}
Also for $x,y\in Q$, define
\[
R(x,y) := R_xR_y R_{xy}\iv 
\qquad \text{and} \qquad
L(x,y) := L_xL_y L_{yx}\iv .
\]
\end{definition}

It is easily seen that $R(x,y) \in \RMlt_1$ and that $\RMlt_1$ is
the group generated
by $\{R(x,y) : x,y \in Q\}$; likewise for the $L(x,y)$ and $\LMlt_1$.

\begin{lemma}
\lemlabel{Mlt}
For any extra loop $Q$:
\begin{itemize}
\item[1.] All permutations in $\RMlt_1$ and $\LMlt_1$ are automorphisms of $Q$.
\item[2.] $R(x,y) R(u,v) = R(u,v) R(x,y)$.
\item[3.] $L(x,y) = R(x,y) = L(y,x) = R(y,x)$
\item[4.] $R(x,y)^2 = I$.
\item[5.] $\RMlt_1 = \LMlt_1$ is a boolean group.
\item[6.] $z R(x,y) = z (x,y,z)$.
\end{itemize}
\end{lemma}

(1) is due to Goodaire and Robinson \cite{GRA}, and (2),(3) are from \cite{KKP};
these are true for all CC-loops. (4) is also from \cite{KKP}, and (5) is immediate
from (2),(3),(4). Also, \cite{KKP} shows that $z L(y,x) = z (x,y,z)\iv$ holds in all CC-loops,
so (6) follows, using (3) and Lemma \lemref{assoc}.

Besides the left and right inner mappings, we have the
middle inner mappings $T_x =  R_x L_x\iv$. In any CC-loop,
the group generated by the middle inner mappings coincides
with the group generated by all inner mappings \cite{DRA1}.

\begin{lemma} 
\lemlabel{T}
In any extra loop $Q$ with $N = N(Q)$ and $A = A(Q)$:
\begin{itemize}
\item[1.] $T_a\in \Aut(Q)$ iff $a\in N(Q)$.
\item[2.] For each $x\in Q$, $\TTT(x) := T_x\res N \in \Aut(N)$.
\item[3.]  $\TTT : Q \to \Aut(N)$ is a homomorphism.
\item[4.] Each $T_x$ maps $A$ onto $A$, so that $A \normal Q$
and $Q/A$ is a group.
\item[5.] Each $(T_x)^2$ is the identity on $A$.
\end{itemize}
\end{lemma}

(1) is from \cite{DRA1}, and holds for all CC-loops.
(2) is due to Goodaire and Robinson \cite{GRA},
and (3) is from \cite{KUNC}.  Both are true for all CC-loops.
$(A)T_x = A$ is due to Fook \cite{FOOK}, and is true for all Moufang loops;
see also Lemma \lemref{assoc-conj} below.
Note that by the remark preceding the lemma, to prove that $A$
is normal, it is sufficient to show that $(A)T_x = A$.
(5) follows from (3) and (4), since $x^2 \in N$, so
$T_{x^2}$ is the identity on $A$ by Lemma \lemref{assoc}.

Our last boolean group is related to two of the others. In an extra
loop $Q$ with $A = A(Q)$, set
\[
A^* := \{ g \in \RMlt : xg \in Ax, \ \forall x \in Q \}
\]
Note that this subgroup of $\RMlt$ is the kernel
of the natural homomorphism $\RMlt(Q) \to \RMlt(Q/A) ;
g \mapsto (Ax \mapsto Axg)$, and so $A^* \normal \RMlt(Q)$.

\begin{lemma}
\lemlabel{A*}
Let $Q$ be an extra loop. Then $A^* = \RMlt_1(Q)\cdot R(A)$,
a direct product. Hence $A^*$ is a boolean group.
\end{lemma}

\begin{proof}
Obviously $R(A)\leq A^*$, and conversely, if $R_a \in A^*$, then
$a\in A$. By Lemma \lemref{Mlt}(6), $\RMlt_1 \leq A^*$. If
$g\in A^*$, write $g = h R_a$ for $h\in \RMlt_1$, $a = 1g$.
Since $h\in A^*$, $R_a\in A^*$, and so $A^* = \RMlt_1 \cdot R(A)$.
Since $A\leq N(Q)$ and $\RMlt_1 \leq \Aut(Q)$, the product
$\RMlt_1 \cdot R(A)$ is direct. Since $A\leq N(Q)$, $R(A)$ is a
boolean group (an isomorphic copy of $A$), and so $A^*$ is a boolean
group by Lemma \lemref{Mlt}(5). 
\end{proof}


\section{The Nucleus}
\seclabel{nuc}

We describe which groups can be nuclei of nonassociative extra loops.

\begin{proposition}
\proplabel{possible-nucs}
For a group $G$, the following are equivalent:
\begin{itemize}
\item[1.] $Z(G)$ contains an element of order 2.
\item[2.] There is a nonassociative extra loop $Q$ with $G = N(Q)$.
\item[3.] There is an extra loop $Q$ with $G = N(Q)$, $|Q : G| = 8$,
and $Z(Q) = Z(G)$.
\end{itemize}
\end{proposition}

\begin{proof}
$(2) \to (1)$ is by Lemma \lemref{assoc}.  Now, assume (1) and
we shall prove (3).
Fix $-1 \in Z(G)$ of order 2, and let 
$C = \{\pm 1, \pm e_1 \dots \pm e_7\}$ be the 16-element Cayley loop.
In the extra loop $G \times C$, let $M = \{(1,1), (-1, -1)\}$.
Note that $M$ is a normal subloop.  Let $Q = (G \times C) / M$.
\end{proof}


\section{Sylow Theorems}
\seclabel{sylow}

We begin by remarking that for extra loops, two possible definitions of
``$p$-loop'' are equivalent. For Moufang loops, the following result is
due to Glauberman and Wright \cite{GLA, GW}. It also holds for
power-associative CC-loops, as follows easily from
(\cite{KKP}, Coro.~3.2, 3.4).

\begin{lemma}
\lemlabel{p-loop}
If $Q$ is a finite extra loop and $p$ is a prime,
then the following are equivalent:
\begin{itemize}
\item[1.] $|Q|$ is a power of $p$.
\item[2.] The order of every element of $Q$ is a power of $p$.
\end{itemize}
\end{lemma}

\begin{definition}
\deflabel{SylowHall}
Let $\pi$ be a set of primes. 
A finite loop $Q$ is a $\pi$-\emph{loop} if the set of prime factors of
$|Q|$ is a subset of $\pi$. If $|Q|$ has prime factorization $|Q| = \Pi_p p^{i_p}$,
then a  \emph{Hall} $\pi$-\emph{subloop} of $Q$ is a subloop of order
$\Pi_{p\in \pi} p^{i_p}$. If $\pi = \{ p\}$, than a Hall $\pi$-subloop
is called a \emph{Sylow} $p$-\emph{subloop}. Let $\Syl_p(Q)$ denote
the set of all Sylow $p$-subloops of $Q$, and let $\Hall_{\pi}(Q)$ denote
the set of all Hall $\pi$-subloops of $Q$.
\end{definition}

Of course, in general, Sylow $p$-subloops and Hall $\pi$-subloops
need not exist. But for extra loops, Sylow $p$-subloops do exist and
satisfy the familiar Sylow Theorems for groups (Theorem \thmref {Sylow} below).
In \S\secref{Hall}, we will show that Hall $\pi$-subloops exist for solvable extra loops
and satisfy P. Hall's Theorem for groups (Theorem \thmref{Hall}).
As a preliminary to both theorems:

\begin{lemma}
\lemlabel{SylowA}
Let $\pi$ be a set of primes with $2\in \pi$, and let $Q$ be a finite extra loop
with $A = A(Q)$.
\begin{enumerate}
\item[1.]  If $P$ is a Hall $\pi$-subloop of $Q$, then $A \le P$.
\item[2.]  If $G$ is a Hall $\pi$-subgroup of $\RMlt(Q)$, then $A^* \le G$.
\end{enumerate}
\end{lemma}

\begin{proof}
Since $A \normal Q$ and is a boolean group, $AP$ is a subloop of $Q$
of order $|A||P|/|A\cap P|$, and so $AP$ is a $\pi$-subloop of $Q$.
By the Lagrange property (Corollary \cororef{Lagrange}),
Hall $\pi$-subloops are maximal $\pi$-subloops,
and so $AP =  P$, establishing (1). The proof for (2) is similar.
\end{proof}

Next we need a minor refinement of the Sylow Theorems for groups.
For a finite group $G$, let $O^p(G)$ denote the subgroup generated
by all elements of order prime to $p$ (\cite{Asch}, p. 5).
Note that $O^p(G) \normal G$.

\begin{lemma}
\lemlabel{SylowB}
Assume that $G$ is a finite group, $p$ is prime, and $P,Q\in \Syl_p(G)$.
Then $Q = x\iv P x$ for some $x\in O^p(G)$.
\end{lemma}

\begin{proof}
If $|G| =  p^m j$, where $p \nmid j$, then
$|O^p(G)| =  p^{\ell} j$, where $0 \le \ell \le m$.
Also $|P\cap O^p(G)| = p^{\ell}$, since
$P\cap O^p(G)\in \Syl_p(O^p(G))$ (\cite{Asch}, (6.4)).
Thus $|P\cdot O^p(G)| = |P| |O^p(G)| / |P\cap O^p(G)| 
= p^m j = |G|$, and so $G = P \cdot O^p(G)$.
Finally,  by the usual Sylow Theorem, let $Q = y\iv P y$,
where $y = ux$, with $u \in P$ and $x \in O^p(G)$.
But then $Q = x\iv P x$.
\end{proof}

\begin{theorem}
\thmlabel{Sylow}
Suppose that $Q$ is a finite extra loop and $|N(Q)| = p^m r$,
where $p$ is prime and $p \nmid r$.  Then
\begin{itemize}
\item[1.]
$|\Syl_p(Q)| = 1 + kp$, where $1 + kp \mid r$.
\item[2.]
If $S$ is a $p$-subloop of $Q$, then there exists $P\in \Syl_p(Q)$
containing $S$.
\item[3.]
If $P_1, P_2\in \Syl_p(Q)$, then
there exists $x\in N(Q)$ such that $P_1 T_x = P_2$,
so that $P_1$ and $P_2$ are isomorphic.
\end{itemize}
\end{theorem}

\begin{proof}
For $p > 2$: By Lemma \lemref{square}(3), every $p$-subloop is
contained in $N$, so the Sylow Theorems for groups can be applied to $N$.

For $p = 2$: The natural homomorphism $[\cdot ]: Q\to Q/A; x\mapsto [x]$
yields a map $[\cdot ]: P\mapsto P/A$ from the set of $2$-subloops $P$ of $Q$
with $A\leq P$ to the set of $2$-subgroups of $Q/A$. 
If $P/A \in \Syl_2(Q/A)$, then $P\in \Syl_2(Q)$, and so 
by Lemma \lemref{SylowA}, $[\cdot ]$ yields a 1-1
correspondence between $\Syl_2(Q)$ and $\Syl_2(Q/A)$.
One can now apply the Sylow Theorems to the group $Q/A$.
To get $x\in N(Q)$ in (3), we apply Lemma \lemref{SylowB}
to $Q/A$ to get $P_1 T_x = P_2$, where $[x] \in O^2(Q/A)$.
Now $x = x_1 \cdots x_n$ where the order of each $[x_i]$,
say $t_i$, is odd.
Then $x_i = a_i z_i$, where $a_i = x_i^{t_i} \in A$ and
$z_i = x_i^{1-t_i} \in N$ since $1-t_i$ is even. Thus each
$x_i\in N$, and so $x\in N$. Finally, that $P_1$ and $P_2$
are isomorphic follows from Lemma \lemref{T}(1).
\end{proof}

Next we relate the Sylow $p$-subloops of an extra loop $Q$ to
the Sylow $p$-subgroups of the right multiplication group $\RMlt(Q)$.

\begin{theorem}
\thmlabel{odd-preSylow}
Let $Q$ be an extra loop with $\RMlt = \RMlt(Q)$.
\begin{enumerate}
\item[1.] If $g\in \RMlt$ has odd order, then $g = R_a$ for some $a\in N(Q)$.
\item[2.] $O^2(\RMlt) \leq R(N(Q))$.
\item[3.] Each subgroup of $\RMlt$ of odd order is isomorphic to a subgroup of $N(Q)$.
\item[4.] $S\mapsto R(S)$ is a 1-1 correspondence between the subloops of
$Q$ of odd order and the subgroups of $\RMlt$ of odd order.
\end{enumerate}
\end{theorem}

\begin{proof}
For $g\in \RMlt$, write (uniquely) $g = h R_a$, where
$a = 1g$ and $h\in \RMlt_1$.
Note that $h R_a h = R_{ah}$ because $h \in \Aut(Q)$ and $h^2 = I$
(Lemma \lemref{Mlt}(1)(5)).
From this plus induction,
$g^{2k} = (R_{ah} R_a)^k$ and 
$g^{2k+1} = h R_a(R_{ah} R_a)^k$ for $k \ge 0$.
Now, the Moufang identity $R_x R_y R_x = R_{xyx}$ plus induction yields
$R_x (R_y R_x)^k = R_{x (yx)^k}$.
Thus, $g^{2k+1} = h R_u$, where $u = a \cdot (ah \cdot a)^k$.
If  $g^{2k+1} = I$ then $h = I$ and $1 = u = a^{2k+1}$, so
$a\in N(Q)$ by Lemma \lemref{square}(4).
This establishes (1), and the rest
follows from (1) and Lemma \lemref{square}(3).
\end{proof}

\begin{theorem}
\thmlabel{even-preSylow}
Let $Q$ be an extra loop. Then $P\mapsto \RMlt_1 \cdot R(P)$
is a 1-1 correspondence between the $2$-subloops of $Q$ containing $A$
and the $2$-subgroups of $\RMlt(Q)$ containing $A^*$.
\end{theorem}

Note that in the theorem, $\RMlt_1 \cdot R(P)$ is not a direct product
of subgroups, but is rather a factorization of a group into a subgroup
and a subset. The multiplication in this group is given by
$h R_a \cdot k R_b = hk R(ak,b) R_{ak\cdot b}$.

\begin{proof}
If $A \leq P\leq Q$, then certainly $A^* \leq \RMlt_1 \cdot R(P)$ by
Lemma \lemref{A*}. Conversely, suppose $G$ is a $2$-subgroup
of $\RMlt$ with $A^* \leq G$,
and set $P = 1G$, the orbit of $G$ through $1\in Q$. Each
$g\in G$ can be uniquely written as $g = h R_a$ for some $h\in \RMlt_1$,
$a = 1g\in P$, and since $\RMlt_1 \leq G$, we have $G = \RMlt_1 \cdot R(P)$.
$|P|$ is a power of $2$, so what remains is to show that $P$ is a subloop.
For $a,b\in P$, $R_a R_b = R(a,b)R_{ab}$, and so
$ab\in P$ as $R(a,b) \leq G$. Similarly, $a\in P$ implies $a\iv\in P$,
which completes the proof.
\end{proof}

\begin{corollary}
\corolabel{Sylow-RMlt}
Let $Q$ be a finite extra loop, and let $p$ be a prime.
Then $\Syl_p(Q)$ is in a 1-1 correspondence with $\Syl_p(\RMlt(Q))$.
\end{corollary}

\begin{proof}
If $p > 2$, then Theorem \thmref{odd-preSylow} yields that $P \mapsto R(P)$ is a
1-1 correspondence between $\Syl_p(Q)$ and $\Syl_p(\RMlt)$.

If $p = 2$, then Theorem \thmref{even-preSylow} and Lemma \lemref{SylowA}(2)
yield that $P \mapsto \RMlt_1\cdot R(P)$ is a 1-1 correspondence
between $\Syl_p(Q)$ and $\Syl_p(\RMlt)$.
\end{proof}


\section{Solvability and Hall $\pi$-subloops}
\seclabel{Hall}

Recall that a loop $Q$ is \emph{solvable} if there exists a normal series
$1 = Q_0 \normal Q_1 \normal \cdots \normal Q_m = Q$ of subloops
$Q_i$ such that each factor $Q_{i+1}/Q_i$ is an abelian group.

\begin{theorem}
\thmlabel{solvable}
An extra loop $Q$ is solvable if and only if $N = N(Q)$ is solvable.
\end{theorem}

\begin{proof}
Since solvability is inherited by subloops, the solvability of $Q$
implies the solvability of $N$. Conversely, if
$1 = N_0 \normal \cdots \normal N_m = N$
is a normal series for $N$, then
$1 = N_0 \normal \cdots \normal N_m \normal Q$
is a normal series for $Q$, since $Q/N$ is an abelian group.
\end{proof}

By Proposition \propref{possible-nucs} and the fact that
the nucleus of a nonassociative extra loop has index at least $8$,
the smallest nonsolvable nonassociative extra loop has order $960$.

\begin{corollary}
\corolabel{Burnside}
Let $Q$ be an extra loop of order $p^a q^b$, where $p,q$ are primes.
Then $Q$ is solvable.
\end{corollary}

\begin{proof}
Since $|N(Q)| = p^c q^d$, the result follows from Burnside's $p^a q^b$-Theorem
for groups (\cite{Asch}, (35.13)) and Theorem \thmref{solvable}.
\end{proof}

This theorem and its corollary actually hold for CC-loops $Q$ because
$Q/N$ is an abelian group by Basarab \cite{BAS} (or see \cite{DRA1, KKP}).
However, the Sylow theorems and P. Hall's Theorem
(cf. \cite{Asch}, (18.5)) can fail in CC-loops,
since the 6-element nonassociative CC-loop does not have a subloop
of order 2.  P. Hall's Theorem for extra loops is:

\begin{theorem}
\thmlabel{Hall}
Let $Q$ be a finite solvable extra loop and $\pi$ a set of primes. Then
\begin{enumerate}
\item $Q$ has a Hall $\pi$-subloop.
\item If $P_1, P_2\in \Hall_{\pi}(Q)$, then there exists $x\in Q$ such
that $P_1 T_x = P_2$.
\item Any $\pi$-subloop of $Q$ is contained in some Hall $\pi$-subloop of $Q$.
\end{enumerate}
\end{theorem}

The proof is similar to that of the Sylow Theorem \thmref{Sylow}.

\begin{proof}
For $2\not\in\pi$: If $S$ is any
$\pi$-subloop of $Q$, then the natural homomorphism $Q \to Q/N$ takes $S$
onto a $\pi$-subloop of a boolean group, so that $S \le N$.  The result then follows
from P. Hall's Theorem applied to the solvable group $N$ (Theorem \thmref{solvable}).

For $2\in \pi$: 
The natural homomorphism $[\cdot ] : Q\to Q/A$ yields a
map $[\cdot ]: P\mapsto P/A$ from the set of $\pi$-subloops
$P$ of $Q$ with $A\leq P$ to the set of $\pi$-subgroups of $Q/A$. 
If $P/A \in \Hall_{\pi}(Q/A)$, then $P\in \Hall_{\pi}(Q)$, and so 
by Lemma \lemref{SylowA}, $[\cdot ]$ restricts to a 1-1
correspondence between $\Hall_{\pi}(Q)$ and $\Hall_{\pi}(Q/A)$.
Now apply P. Hall's Theorem to the solvable group $Q/A$.
\end{proof}


\section{The Center}
\seclabel{center}

\begin{theorem}
\thmlabel{center}
If $Q$ is a nonassociative extra loop and $A(Q)$ is finite, then
$|Z(Q)\cap A(Q)| > 1$.
\end{theorem}

\begin{proof}
Applying Lemma \lemref{T}, define $\TTT': Q \to \Aut(A)$ by
$\TTT'(x) = T_x\res A$. By Lemma \lemref{assoc}, $\TTT'(x) = I$ for $x \in N$.
Thus, via $\TTT'$, the boolean group $Q/N$ acts on the boolean group $A$.
Since $|A|$ is even and the size of each orbit is a power of $2$,
there must be some $a \in A \backslash \{1\}$ which is fixed by
this action.  Then $a\in Z(Q)$.
\end{proof}

This can fail when $A(Q)$ is infinite; see Example \exref{inf}.

\begin{lemma}
\lemlabel{assoc-conj}
In an extra loop,
\[
(x,y,zt) =(x,y,tz) = (x,y,z) \cdot (x,y,t)T_z = (x,y,z)T_t \cdot (x,y,t) \ \ .
\]
\end{lemma}

\begin{proof}
Applying Lemma \lemref{Mlt}, we have
$z R(x,y) = z (x,y,z)$,
$t R(x,y) = t (x,y,t)$, and
$ z R(x,y) \cdot t R(x,y)  = (zt) R(x,y) = zt \cdot (x,y,zt)$,
so 
\[
z (x,y,z) \cdot t (x,y,t)  = zt \cdot (x,y,zt) \ \ .
\]
Since associators are in the nucleus,
we get $(x,y,z)T_t \cdot  (x,y,t) =  (x,y,zt)$.
Also, $(x,y,tz) = (x,y,zt)$ by Lemma \lemref{assoc},
since $Q/N$ is abelian, so $tz  \in Nzt$.
\end{proof}

Since $(x,y,t)T_z = (x,y,z) \cdot (x,y,zt) $, we have,
in the case of extra loops, another proof
of Fook's result (Lemma \lemref{T}.3) that $(A)T_z = A$.
Lemma \lemref{assoc-conj} yields:

\begin{lemma}
\lemlabel{assoc-commutes}
In an extra loop, $z$ commutes with $(x,y,t)$ iff
$t$ commutes with $(x,y,z)$ iff
$(x,y,z) (x,y,t) = (x,y,zt)$.
\end{lemma}

\begin{lemma}
\lemlabel{assoc-nuc}
If $Q$ is an extra loop, with $a = (x,y,z)$, then
$a \in Z(\sbl{\{x,y,z\} \cup N})$, and
$A(\sbl{\{x,y,z\} \cup N}) = \{1,a\}$
\end{lemma}

\begin{proof}
$a \in N$ implies that $T_a$ is an automorphism of $Q$
(Lemma \lemref{T}), so that $\{s\in Q : sa = as\}$ is a subloop
of $Q$, and this subloop contains all elements of $\{x,y,z\} \cup N$
by Lemma \lemref{assoc}, which also implies that
$(u,v,w) \in \{1,a\}$ for all $u,v,w \in \{x,y,z\} \cup N$.
Then $A(\sbl{\{x,y,z\} \cup N}) \subseteq \{1,a\}$ follows
by using Lemma \lemref{assoc-conj}.
\end{proof}

\begin{lemma}
\lemlabel{assoc-index}
If $Q$ is an extra loop, then $|A(Q) : A(Q) \cap Z(Q)| \notin \{2,4,8\}$.
\end{lemma}

\begin{proof}
Set  $Z = A(Q) \cap Z(Q)$, and define
$\TTT': Q \to \Aut(A)$, as in the proof of Theorem \thmref{center}.
Assume that $|A : Z|  > 1$.
Fix $e_1,e_2,e_3 \in Q$ with $(e_1,e_2,e_3) \notin Z$, and then fix
$e_4\in Q$ such that
$(e_1,e_2,e_3) \TTT'(e_4) \neq (e_1,e_2,e_3)$.
Define
\[
q_1 := (e_2,e_3,e_4) \quad
q_2 := (e_1,e_3,e_4) \quad
q_3 := (e_1,e_2,e_4) \quad
q_4 := (e_1,e_2,e_3) \ \ .
\]
By Lemmas \lemref{assoc-commutes} and \lemref{assoc},
$q_i \TTT'(e_j) = q_i$ iff $j \ne i$.
Now, let $q_S = \prod_{i\in S} q_i$ for $S \subseteq \{1,2,3,4\}$,
and observe that $q_S \TTT'(e_j) = q_S$ iff $j \notin S$,
so that the $q_S$ are all in distinct cosets of $Z$.
Thus, $|A:Z| \geq  16$.
\end{proof}

\begin{theorem}
\thmlabel{assoc-512}
If $Q$ is a finite extra loop with some associator not contained in $Z(Q)$,
then $|A(Q)| \ge 32$ and $|Q : N(Q)| \ge 16$, so that $512 \mid |Q|$.
\end{theorem}
\begin{proof}
$|Q : N| \ge 16$ follows from Lemma \lemref{assoc-nuc}.
$| A(Q) \cap Z(Q)| \ge 2$ follows from Theorem \thmref{center},
so $|A(Q)| \ge 32$  follows from Lemma \lemref{assoc-index},
so that $512 \mid |Q|$.
\end{proof}

The ``$512$'' is best possible; see Example \exref{512}.
The construction there
is suggested by the proof of Lemma \lemref{assoc-index}.
We shall get $A(Q) = N(Q) = \sbl{q_0, q_1, q_2, q_3, q_4}$, of order $32$,
$Q/N = \sbl{[e_1],[e_2],[e_3],[e_4]}$, of order $16$,
and $Z(Q) = \{1, q_0\}$.


\section{Extension}
\seclabel{ext}

Say we are given an abelian group $(G, +)$ and a boolean group
$(B, +)$, and we wish to construct all extra loops $Q$ such that
$G \normal Q$, $G \le N(Q)$, and  $Q/G \cong B$.
We may view this as an extension problem; see
\cite{CPS} \S II.3, p.~35.

Assuming that we already have $Q$, let $\pi : Q \to B$ be the natural
quotient map.  By the Axiom of Choice, we can assume that $B$
is a section; that is, $B$ is a subset of $Q$ and
$\pi\res B$ is the identity function.
Then for $a,b \in B$, we have the loop product $a \cdot b$ from 
$Q$ and the abelian group sum $a + b \in B$.
Since $a \cdot b$ and $a + b$ are in the same left coset of $G$,
there is a function $\psi : B \times B \to G$ with
$a \cdot b = (a + b) \psi(a,b)$.
We may assume that the identity element of $B$ is the $1$ of $Q$, so that
$\psi(1,a) = \psi(a,1) = 1$.
Each $T_a \res G \in \Aut(G)$.
Also, the map $x \mapsto T_x \res G$ is a homomorphism from $Q$
to $\Aut(G)$, and is the identity map on $G$ (since $G$ is abelian),
so it defines a homomorphism: $B \to \Aut(G)$.
Every element of $Q$ is
in some left coset of  $G$, so it can be expressed uniquely in the
form $a u$, with $a \in B$ and $u \in G$.
Since  $G \le N(Q)$, we can compute the product of two elements of
this form as
$ a u \cdot  b v =  a  b \cdot u T_{ b} v =
(a + b) \cdot  \psi(a,b)  ( u T_{ b}) v $.
In particular, for $b \in B$, 
$b^2 =  b \cdot b = (b + b) \cdot  \psi(b,b) = \psi(b,b)$.

Turning this around, and converting to additive notation,

\begin{definition}
\deflabel{extension}
Suppose we are given:
\begin{itemize}
\item [1.]
An abelian group $(G, +)$ and a boolean group $(B, +)$.
\item [2.]
A map $\psi : B \times B \to G$  with $\psi(0,a) = \psi(a,0) = 0$.
\item [3.]
A homomorphism, $a \mapsto \tau_a$, from $B$ to $\Aut(G)$.
\end{itemize}
Then $B\ltimes_\tau^\psi G$ denotes the set $B \times G$ given
the product operation:
\[
(a, u) \cdot (b,v) = (a + b, \;  \psi(a,b)  +   u \tau_{b} + v ) \ \ .
\]
$B \ltimes_\tau G$ denotes $B\ltimes_\tau^\psi G$ 
in the case that $\psi(a,b) = 0$ for all $a,b$.
\end{definition}

Then $B \ltimes_\tau G$ is a group, and is the usual semidirect product.

\begin{lemma}
\lemlabel{extension-loop}
$B \ltimes_\tau^\psi G$ is always a loop with identity element $(0,0)$.
The map $u \mapsto (0,u)$ is an isomorphism from $G$
onto $\{0\}\times G \normal B \ltimes_\tau^\psi G$.
\end{lemma}

\begin{proof}
We can solve the equations
$(a, u) \cdot (b,v) = (c,w)$ for $(b,v)$ or $(a,u)$:
\begin{align*}
(a,u) \ld (c,w) &= (a + c, \; w - \psi(a,a+c) - u \tau_a\tau_c  ) \\
(c,w) \rd (b,v) &= (b + c, \; w\tau_b - \psi(b+c,b)\tau_b - v \tau_b )  \ \ .
\end{align*}
Here, we have simplified the expression using the facts that $B$ is
boolean and the map $b \mapsto \tau_b$ is a homomorphism.
This proves that $B \ltimes_\tau^\psi G$ is a loop.  
$\{0\}\times G$ is a normal subloop because the map
$(a,u) \mapsto a$ is a homomorphism.
\end{proof}

It is fairly easy to calculate, in terms of $\psi$ and $\tau$,
what is required for $B \ltimes_\tau^\psi G$ to satisfy various
properties, such as the inverse property, the Moufang law, etc.
In the case of extra loops, we shall use the conditions of
Lemma \lemref{extra-assoc} on the associators;
some of these conditions can be verified immediately:

\begin{lemma}
\lemlabel{assocpsi}
Let $Q = B \ltimes_\tau^\psi G$.  Then $A(Q) \le \{0\}\times G \le N(Q)$.
\end{lemma}

\begin{proof}
To compute the associators, we solve:
\[
[(a,u) \cdot (b,v) (c,w)] \;\cdot\; \big((a,u),(b,v),(c,w)\big)
\ =\  (a,u)(b,v) \cdot (c,w) \ \ .
\]
First, we compute both associations:
\begin{align*}
(a,u) \cdot (b,v) (c,w) &= (a,u) (b + c,\; \psi(b,c) + v\tau_c + w)   \\
     &= (a+b+c,\; \psi(a, b+c) + u \tau_b \tau_c +  \psi(b,c) + v\tau_c + w) \\
(a,u)(b,v) \cdot (c,w) &= (a+b,\; \psi(a,b) + u \tau_b + v) \cdot (c,w) \\
      &= (a+b+c,\; \psi(a+b,c) +
                \psi(a,b)\tau_c + u \tau_b\tau_c + v\tau_c + w) \ \ .
\end{align*}
So,
\[
\big((a,u),(b,v),(c,w)\big) 
\ = \   \big(0, \;
\psi(a+b,c) + \psi(a,b)\tau_c  - \psi(a, b+c) -  \psi(b,c) \big) \ \ .
\]
Observe that this depends only on $a,b,c$, and has value $0$ 
if any of $a,b,c$ are $0$, so that
$\{0\}\times G \le N(Q)$, and all $(x,y,z) \in \{0\}\times G$.
\end{proof}

We now consider in more detail the case when both $B$ and $G$ are
boolean.  We shall in fact start with $\tau$ and the desired
associator map $\alpha : B^3 \to G$, where
$\big(0, \alpha(a,b,c) \big)$ denotes the intended value of
$\big((a,u),(b,v),(c,w)\big)$ for some (any) $u,v,w \in G$.
We plan to construct $\psi$ from $\alpha$ and $\tau$.
This is useful because $\alpha$ is determined
by its values on a basis for $B$.  We need to assume some
conditions on $\alpha$ suggested by Lemmas \lemref{assoc-conj}
and \lemref{assoc}:

\begin{lemma}
\lemlabel{extendalpha}
Suppose that $G$ and $B$ are boolean groups and $\tau \in \Hom(B, \Aut(G))$.
Let $E$ be a basis for $B$, and assume that
$\alpha: E^3 \to G$ satisfies the equations:
\begin{itemize}
\item[H1.] $(\alpha(a_1,b,c)) \tau_{a_2} +     \alpha(a_2,b,c)  =
      \alpha(a_1,b,c)    + (\alpha(a_2,b,c)) \tau_{a_1} $.
\item[H2.] $(\alpha(a,b_1,c)) \tau_{b_2} +     \alpha(a,b_2,c)  =
      \alpha(a,b_1,c)    + (\alpha(a,b_2,c)) \tau_{b_1} $.
\item[H3.] $(\alpha(a,b,c_1)) \tau_{c_2} +     \alpha(a,b,c_2)  =
      \alpha(a,b,c_1)    + (\alpha(a,b,c_2)) \tau_{c_1} $.
\item[F1.] $(\alpha(a,b,c)) \tau_{a}  = \alpha(a,b,c) $.
\item[F2.] $(\alpha(a,b,c)) \tau_{b}  = \alpha(a,b,c) $.
\item[F3.] $(\alpha(a,b,c)) \tau_{c}  = \alpha(a,b,c) $.
\end{itemize}
Then $\alpha$ extends uniquely to a map $\oa : B^3 \to G$ satisfying
these same equations for all elements of $B$, together with
\begin{itemize}
\item[P1.] $\oa(a_1+a_2,b, c) =
(\oa(a_1,b,c)) \tau_{a_2} +  \oa(a_2,b,c)$.
\item[P2.] $\oa(a,b_1 + b_2, c) =
(\oa(a,b_1,c)) \tau_{b_2} +  \oa(a,b_2,c)$.
\item[P3.] $\oa(a,b, c_1 + c_2) =
(\oa(a,b,c_1)) \tau_{c_2} +  \oa(a,b,c_2)$.
\end{itemize}
If $\alpha$ is symmetric, then the same holds for $\oa$.
If in addition,
$\alpha$ satisfies $\alpha(a,a,b) = 0$ for all $a,b\in E$,
then $\oa(a,a,b) = 0$ for all $a,b\in B$.
\end{lemma}

\begin{proof}
First, fix $a,b \in E$, and consider the map
$\varphi:  E \to B \ltimes_\tau G$ 
defined by $\varphi(c) = (c, \alpha(a,b,c))$.
H3 says that $\varphi(c_1) \varphi(c_2) = \varphi(c_2) \varphi(c_1)$,
and F3 says that each $(\varphi(c))^2 = 1$.
It follows that $\varphi$ extends uniquely to a homomorphism
$\varphi' : B \to B \ltimes_\tau G$; then
$\varphi'(c) = (c, \alpha'(a,b,c))$.

Doing this for every $a,b \in E$, we get 
$\alpha' : E \times E \times B \to G$, which is the unique extension
of $\alpha$ satisfying H3,F3,P3.
But then it is easily seen that $\alpha'$ satisfies H1,H2,F1,F2 also.
$\alpha'$ is computed inductively using P3;
the purpose of $\varphi$ was just to prove that this computation
yields a well-defined function.

Repeating this on the second coordinate yields
$\alpha'' : E \times B \times B \to G$, which is the unique extension
of $\alpha$ satisfying H2,H3,F2,F3,P2,P3.  Doing it again yields $\oa$.

If $\alpha$ is symmetric, then the symmetry of $\oa$ follows from the
uniqueness of $\oa$.  Finally, assume in addition that
$\alpha(a,a,b) = 0$ holds on $E$.  First, for each $e \in E$,
note that $\{b \in B : \oa(e,e,b) = 0\}$ is a subgroup of $B$,
so that $\oa(e,e,b) = 0$ for all $b \in B$.
Then, for each fixed $b \in B$, $\{a \in B : \oa(a,a,b) = 0\}$ is a
also a subgroup, so that $\oa(a,a,b)$ for all $a,b \in B$.
\end{proof}

We now analyze the special case that in 
$Q = B \ltimes_\tau^\psi G$,
the elements of $E \times \{0\}$ all have order $2$ 
and all commute with each other.
We can then use $\alpha$ to compute the correct $\psi$.
Observe first:

\begin{lemma}
\lemlabel{right-assoc}
In an extra loop $Q$, suppose that the elements $x_1, x_2, \ldots , x_n$
all pairwise commute.  Let $\pi$ be a permutation of $\{1,2, \ldots , n\}$.
Then $x_1 \cdot x_2 \cdot \cdots  \cdot x_n = 
x_{\pi(1)} \cdot x_{\pi(2)} \cdot \cdots  \cdot x_{\pi(n)}$,
where both products are right-associated.
\end{lemma}

\begin{proof}
It is sufficient to prove  $x \cdot y z = y \cdot xz$ when $xy = yx$,
and this follows by $x \cdot y z = x y \cdot z \cdot (x,y,z) = 
 y x \cdot z \cdot (x,y,z) = y \cdot xz$.
\end{proof}

Thus, if the elements of $E \times \{0\}$ all commute, then the value
of a right-associated product from $E \times \{0\}$ must be independent of
the order in which that product is taken.  This will simplify
the form of $\psi$.  If the elements of $E \times \{0\}$ also have order $2$
in $Q$, then it is easy to say what properties $\alpha$ must
satisfy:

\begin{theorem}
\thmlabel{build-extra}
Suppose that we are given boolean groups $G$ and $B$, with $E \subset B$
a basis for $B$.  Suppose that we also have $\tau \in \Hom(B, \Aut(G))$
and a map $\alpha: E^3 \to G$ satisfying:
\begin{itemize}
\item[1.] $\alpha$ is invariant under permutations of its arguments.
\item[2.] $\alpha(e_1, e_1, e_2) = 0$.
\item[3.] $(\alpha(e_1,e_2,e_3)) \tau_{e_4} +     \alpha(e_1,e_2,e_4)  =
      \alpha(e_1,e_2,e_3)    + (\alpha(e_1,e_2,e_4)) \tau_{e_3} $.
\end{itemize}
Then there is a unique $\psi:  B \times B \to G$ satisfying:
\begin{itemize}
\item[a.]  $\psi(0,a) = \psi(a,0) = 0$ for all $a\in B$.
\item[b.] $Q := B \ltimes_\tau^\psi G$ is an extra loop.
\item[c.] In $Q$, whenever $e_1,e_2,e_3 \in E$, we have
$(e_1, 0) \cdot (e_1, 0) = 0$,
$(e_1, 0) \cdot (e_2, 0) = (e_2, 0) \cdot (e_1, 0)$,
and the associator
$\big( (e_1, 0) , (e_2, 0) ,  (e_3, 0)\big) = \big(0, \alpha(e_1,e_2,e_3)\big)$.
\item[d.] $\psi(e,b ) = 0$  whenever $e\in E$.
\end{itemize}
\end{theorem}

Condition (d) expresses the intent that the
elements of the section be right-associated products from $E$.

\begin{proof}
Note that (1 -- 3) implies that
$(\alpha(e_1,e_2,e_3)) \tau_{e_1}  = \alpha(e_1,e_2,e_3)$.

By Lemma \lemref{extendalpha}, $\alpha$ extends uniquely to a symmetric
map $\oa: B^3 \to G$ satisfying the conditions H$i$,F$i$,P$i$ there.
For the uniqueness part of 
the theorem, we note that assuming that
$B \ltimes_\tau^\psi G$ is an extra loop, 
this $\oa$ must indeed yield
the associator; that is, by condition (c) and Lemma
\lemref{assoc-conj}, we have:
\[
\big((a,u),(b,v),(c,w)\big)  = \big(0, \oa(a,b,c)\big) \ \ .
\]
Then, by the computation in
the proof of Lemma \lemref{assocpsi}, we get:
\[
\oa(a,b,c) =
 \psi(a + b,c) + \psi(a,b)\tau_c + \psi(a, b+c) + \psi(b,c)   \ \ .
\]
Consider the case where $a = e \in E$.
Then condition (d) implies that $\psi(e,b) = \psi(e, b+c) = 0$, so we get
$\psi(e+b, c) = \psi(b,c) + \oa(e,b,c)$.
Repeating this, we see that  for $e_1, \ldots, e_n \in E$,
\[
\psi(e_1 + \cdots  + e_n , c)
\ =\  \sum_{j=1}^n \oa\big(e_j,\; \sum_{k<j}e_k,\; c\big) \ \ .
\eqno{(*)}
\]
For example, 
\begin{align*}
\psi(e_1 + e_2, c) = &\ \oa(e_2, e_1, c) \\
\psi(e_1 + e_2 + e_3, c) = &\ \oa(e_2, e_1, c) + \oa(e_3, e_1+e_2, c) = \\
&\ \oa(e_2, e_1, c) + (\oa(e_3, e_1, c))\tau_{e_2} +  \oa(e_3, e_2, c)
\end{align*}
This proves the uniqueness of $\psi$.  To prove existence,
one can take $(*)$ as a definition of $\psi$ (after proving
that it is well-defined), and then prove that it yields an
extra loop with the correct associators.

To prove that it is well-defined, fix $c$ and define,
$\Psi_n = \Psi_n^{(c)}: E^n \to B$ for $n\ge 1$ so that
\begin{align*}
\Psi_1(e) &= 0. \\
\Psi_{n+1}(e_0, e_1, \ldots, e_n) &=  \Psi_{n}(e_1, \ldots, e_n) +
 \oa(e_0, e_1 + \cdots + e_n, c) \ \ .
\end{align*}
It is easy to see that
$\Psi_2(e,e) = 0$ and that
$\Psi_{n+2}(e,e, e_1, \ldots, e_n) = \Psi_{n}(e_1, \ldots, e_n)$.
We need to prove that each $\Psi_n$ is invariant under permutations
of its arguments.  Then, it will be unambiguous to define
$\psi( e_1 + \cdots + e_n, c)= \Psi_n^{(c)}( e_1, \ldots, e_n)$.
To prove invariance under permutations, we induct on $n$;
for the induction step, it is sufficient to prove that
$\Psi_{n+2}(e,e', e_1, \ldots, e_n) = \Psi_{n+2}(e',e, e_1, \ldots, e_n)$,
and this follows from the fact that
\begin{align*}
\oa(e, e' + b, c) + \oa(e', b, c) &=
(\oa(e, e', c))\tau_b + \oa(e, b, c) + \oa(e', b, c)  \\
&= \oa(e', e + b, c) + \oa(e, b, c) \ \ .
\end{align*}

Now that we have $\psi$ defined, we need to check that our given
$\oa(a,b,c)$ is really the true associator.  Use
$\big(0, (a,b,c)\big)$ to denote
$\big((a,u),(b,v),(c,w)\big)$ for some (any) $u,v,w \in G$;
then, as in the
proof of Lemma \lemref{assocpsi},
\[
(a,b,c) =
 \psi(a + b,c) + \psi(a,b)\tau_c + \psi(a, b+c) + \psi(b,c)   \ \ .
\]
We prove $\oa(a,b,c) = (a,b,c)$ by induction on the number of basis
elements needed to add up to $a$.  If $a = 0$, then $\oa(a,b,c) = (a,b,c) = 0$.
For the induction step, note that
$\oa(e+a, b, c) - \oa(a,b,c) = \oa(e,b,c) \tau_a$, which is the same
as $(e+a, b, c) - (a,b,c)$, since using 
$\psi(e+b, c) = \psi(b,c) + \oa(e,b,c)$, we get:
\begin{align*}
&(e+a, b, c) - (a,b,c) = \oa(e, a+b, c) + \oa(e,a,b)\tau_c + \oa(e,a,b+c)=\\
&\qquad \oa(e,b,c)\tau_a + \oa(e,a,c) +  \oa(e,a,b)\tau_c +
          \oa(e,a,b)\tau_c + \oa(e,a,c) = \oa(e,b,c)\tau_a \  .
\end{align*}

Now that we have identified $\oa(a,b,c)$ as the associator, it
is easy to prove that $Q$ is an extra loop by verifying the conditions
in Lemma \lemref{extra-assoc}.
(2) and (3) are clear from Lemma \lemref{assocpsi}.
(1) ($Q$ is flexible) holds because $\oa(a,b,a) = 0$,
and (4) holds because $\oa$ is symmetric.
For (5), we must check that $\big(0, \oa(a,b,c)\big)$ commutes with $(a,u)$,
and this follows from the fact that $(\oa(a,b,c))\tau_a = \oa(a,b,c)$.
\end{proof}

We now describe three examples.

If $|G| = 2$ and $|B| = 8$
(so $E = \{e_1,e_2,e_3\}$), there is only 
one non-associative option.
$\alpha(e_1,e_2,e_3)$ must be the non-identity element of $G$,
and each $\tau_x$ must be $I$.
This extra loop of order $16$ is the opposite extreme from the Cayley loop
(where the elements outside the nucleus have order $4$ and
anticommute).

\begin{example}
\exlabel{512}
There is an extra loop $Q$ of order $512$
such that $Q/Z(Q)$ is nonassociative.
\end{example}

\begin{proof}
Let $E = \{e_1, e_2, e_3, e_4\}$ and
$G = \sbl{q_0, q_1, q_2, q_3, q_4}$, so that $|Q| = 512$.
Define $\tau$ so that $q_0 \tau_{e_k} = q_0$ and
$q_j \tau_{e_k} = q_j + \delta_{j,k} q_0$ for $j,k \in \{1,2,3,4\}$;
then $Z(Q)$ will be $\{ (0,0),  (q_0, 0) \}$.
Define $\alpha$ so that $\alpha(e_i,e_j,e_k) = q_\ell$ whenever
$i,j,k,\ell \in \{1,2,3,4\}$ are distinct.
\end{proof}

The $\psi$ of this
example was first found using McCune's program Mace4 \cite{McM},
and the abstract discussion of this section was then obtained
by reverse engineering.

\begin{example}
\exlabel{inf}
There is an infinite nonassociative extra loop $Q$
with $Z(Q) = \{1\}$.
\end{example}

\begin{proof}
Let $B$ be any infinite boolean group, and
we use a wreath product construction.
$B$ acts on $(\ZZ_2)^B$ by permuting the indices; that is,
for $u : B \to \ZZ_2$, let $((u)\tau_a)(b) = u(a+b)$.
Let $G = \{u \in (\ZZ_2)^B : |u\iv\{1\}| < \aleph_0\}$; so $G$
is a direct sum of $|B|$ copies of $\ZZ_2$ (and is hence isomorphic
to $B$, since $\dim(B) = |B|$).
Since $B$ is infinite, $B \ltimes_\tau G$ (and hence
also $B \ltimes_\tau^\psi G$) will have trivial center.

Let $E$ be a basis for $B$.  For $e_1,e_2,e_3 \in E$, let
$\alpha(e_1,e_2,e_3) = 0$  unless $e_1,e_2,e_3$ are distinct,
in which case $\alpha(e_1,e_2,e_3)$ is the element of 
$G \le (\ZZ_2)^B$ which is $1$ on the $8$ members 
of $\sbl{e_1,e_2,e_3}$ and $0$ elsewhere.
To verify condition (3), we let
$u = (\alpha(e_1,e_2,e_3)) \tau_{e_4} +  \alpha(e_1,e_2,e_4)$
and let $v = \alpha(e_1,e_2,e_3)    + (\alpha(e_1,e_2,e_4)) \tau_{e_3} $,
and consider cases:
If $e_1 = e_2$, then $u = v = 0$, so assume that $e_1 \ne e_2$.
If $e_3 \in \{e_1,e_2\}$, then $u = v = \alpha(e_1,e_2,e_4)$, and
if $e_4 \in \{e_1,e_2\}$, then $u = v = \alpha(e_1,e_2,e_3)$, 
so assume also that $\{e_3,e_4\} \cap \{e_1,e_2\} = \emptyset$.
If $e_3 = e_4$ then $u = v = 0$.
In the remaining case,
$e_1, e_2, e_3, e_4$ are all distinct; then
both $u,v$ are $1$ on the $16$ members
of $\sbl{e_1,e_2,e_3,e_4}$ and $0$ elsewhere.
\end{proof}

\section{Semidirect Products}
\seclabel{semi}

The loop $B\ltimes_\tau^\psi G$ from Definition \defref{extension}
is not really a semidirect product, since it need not contain
an isomorphic copy of $B$.  
If we delete the $\psi$, we get a true semidirect product.
Following Robinson \cite{ROB}:

\begin{definition} 
\deflabel{semider}
Let $B, G$ be loops, and assume that $\tau \in \Hom(B, \Aut(G))$.
Then $B \ltimes_\tau G$ denotes the set $B\times G$ given the product operation:
\[
(a, u) \cdot (b, v) = (a b,\, (u)\tau_b  \cdot v) \ \ .
\]
We write $B \ltimes G$ when $\tau$ is clear from context.
\end{definition}

It is easily verified that $B \ltimes G$ is a loop,
with identity element $(1,1)$,
but $B \ltimes G$ need not inherit all the properties satisfied by $B$ and $G$.
The general situation for extra loops was discussed in \cite{ROB}.  Here,
we consider only an easy special case:

\begin{lemma}
\lemlabel{inherit}
Assume that $\tau \in \Hom(B, \Aut(G))$,
$B$ is an extra loop, and $G$ is a group.
Then $B \ltimes_\tau G$ is an extra
loop, and the inverse is given by
$(a,u)\iv = (a\iv, (u\iv)\tau_{a\iv})$.
\end{lemma}
\begin{proof}
Note that $(a,u) \cdot (a\iv, (u\iv)\tau_{a\iv}) = (1,1)$.
We verify the extra loop equation $(xy\cdot z) \cdot x = x \cdot (y\cdot zx)$,
setting $x = (a,u)$, $y = (b,v)$, $z = (c,w)$:
\begin{align*}
((a,u) (b,v) \cdot (c,w)) \cdot (a,u) &=
\big( (ab\cdot c) \cdot a  , 
    \  (u)\tau_{bca} \cdot (v)\tau_{ca} \cdot(w) \tau_a \cdot u \big) \\
(a,u) \cdot ((b,v)\cdot (c,w) (a,u))  &=
\big( a \cdot (b\cdot ca) ,
    \  (u)\tau_{bca} \cdot (v)\tau_{ca} \cdot(w) \tau_a \cdot u \big) \\
\end{align*}
These are clearly equal, since $B$ is an extra loop.
In writing these equations, we used the facts that $G$ is associative,
and that $\Aut(G)$ is associative and $\tau$ is a homomorphism,
so that the notation $\tau_{bca}$ is unambiguous,
even though $b\cdot ca$ need not equal $bc \cdot a$.
\end{proof}

Of course,
the same reasoning will work for other equations
which are weakenings of the associative law;
for example, if $B$ is Moufang and $G$ is a group, then 
$B \ltimes G$ is Moufang.

In some cases, we can prove that every extra loop of a given
order is a semidirect product:

\begin{lemma}
\lemlabel{issemi}
Suppose that $Q$ is a finite extra loop and $N = N(Q)$ is abelian.
Then $Q$ is isomorphic to $B \ltimes_\tau G$, where
$B \in \Syl_2(Q)$, $G = O^2(N)$, $\tau_a = T_a \res G$,
and each $(\tau_a)^2 = I$.
\end{lemma}
\begin{proof}
Say $|Q| = 2^n r$, where $r$ is odd, so $|B| = 2^n$.
Then $|N| = 2^m r$ for some 
$m \le n$, and $|B \cap N| = 2^m$.  Since $N$ is abelian,
it is an internal direct sum of $B \cap N$ and $G = O^2(N)$,
which must have order $r$.  Then $Q = BG$, since $B \cap G = \{1\}$.
Furthermore, each $T_a$ maps $G$ to  $G$ because
$T_a \in \Aut(N)$ and $G$ is a characteristic subgroup of $N$.
Then $Q \cong B \ltimes_\tau G$ follows.  Also, 
$(\tau_a)^2 = \tau_{a^2} = I$ because $a^2 \in N$, which is abelian.
\end{proof}

\begin{lemma}
\lemlabel{sixteenp}
Suppose that $Q$ is a nonassociative
extra loop of order $16p$, where $p$ is an
odd prime.  Then $N(Q) \cong \ZZ_2 \times \ZZ_p$.
\end{lemma}
\begin{proof}
$|Q : N| \ge 8$ because any $\sbl{\{x,y\} \cup N}$ is associative,
and $Z(N)$ contains an element of order $2$ by Lemma \lemref{assoc},
so $|N| = 2p$ and $N$ cannot be the dihedral group, so $N$ must be
$\ZZ_2 \times \ZZ_p$.
\end{proof}

Combining Lemmas \lemref{issemi} and \lemref{sixteenp},
we see that such $Q$ must be of the form $B \ltimes_\tau \ZZ_p$,
where $B$ is one of the five extra loops of order $16$
and each $\tau_a \in \{1,-1\} \le \Aut(\ZZ_p)$;
this is because $(\tau_a)^2 = I$, and the only element
of $\Aut(\ZZ_p)$ of order $2$ is the map $u \mapsto -u$.
We shall now show that the number of such loops is independent of $p$.
Obviously, $\Hom(B, \{1,-1\})$ does not depend on $p$, but
different homomorphisms can result in isomorphic loops, so we 
must show that for $\tau, \sigma \in \Hom(B, \{1,-1\})$,
the question of whether $B \ltimes_\tau \ZZ_p \cong B \ltimes_\sigma \ZZ_p$
does not depend on $p$:

\begin{lemma}
\lemlabel{allsixteenp}
If $B$ is a finite extra $2$-loop and $\tau, \sigma \in \Hom(B, \{1,-1\})$,
say $\tau \sim \sigma$ iff there is an $\alpha \in \Aut(B)$
with $\tau = \alpha \sigma$.  Let $p$ be an odd prime.
Then, identifying $\{1,-1\} \le \Aut(\ZZ_p)$,
$B \ltimes_\tau \ZZ_p \cong B \ltimes_\sigma \ZZ_p$
iff  $\tau \sim \sigma$.
\end{lemma}
\begin{proof}
If $\tau = \alpha \sigma$, then define
$\Phi: B \ltimes_\tau \ZZ_p \to  B \ltimes_\sigma \ZZ_p$ by
$(a,u)\Phi = ((a)\alpha, u)$.
To verify that $\Phi$ is an isomorphism, use
\begin{align*}
( (a,u) \cdot_\tau (b,v) ) \Phi \ &=\  (ab,\  (u)\tau_b + v) \Phi
\ = \ ((ab)\alpha,\  (u)\tau_b + v) \\
(a,u)\Phi \cdot_\sigma (b,v)\Phi\  &= \
((a)\alpha, u) \cdot_\sigma ((b)\alpha, v)
\ =\ ((a)\alpha \cdot (b)\alpha,\  (u)\sigma_{(b)\alpha} + v) \ \ ,
\end{align*}
and these are equal because $\tau_b$ (i.e., $(b)\tau$)
is the same as $\sigma_{(b)\alpha}$ (i.e., $(b)\alpha\sigma$).

Conversely, suppose we are given an isomorphism 
$\Phi: B \ltimes_\tau \ZZ_p \to  B \ltimes_\sigma \ZZ_p$.
Then $\Phi( B \times \{0\}) \in \Syl_2( B \ltimes_\sigma \ZZ_p)$.
But also $(B \times \{0\}) \in \Syl_2( B \ltimes_\sigma \ZZ_p)$,
and $\Aut( B \ltimes_\sigma \ZZ_p)$ acts transitively on 
the set of Sylow 2-subloops by Theorem \thmref{Sylow}.
Thus, composing $\Phi$ with an automorphism, we may assume
WLOG that $\Phi( B \times \{0\}) = B \times \{0\}$.
Also, $\Phi( \{1\} \times \ZZ_p) =  \{1\} \times \ZZ_p$ because
$\{1\} \times \ZZ_p$ is the only subloop of $B \ltimes_\sigma \ZZ_p$
isomorphic to $\ZZ_p$.  So, we have
$(a,0)\Phi = ((a)\alpha, 0)$ and $(1,u)\Phi = (1, (u)\beta)$
for some $\alpha \in \Aut(B)$ and $\beta \in \Aut(\ZZ_p)$.
Since $(a,u) = (a,0) \cdot (1,u)$, we also have
$(a,u)\Phi = ((a)\alpha, (u)\beta)$.  Furthermore, the map
$(c,w) \mapsto (c, (w)\beta\iv)$ is an automorphism
of $B \ltimes_\sigma \ZZ_p$, since $\Aut(\ZZ_p) \cong \ZZ_{p-1}$ is
abelian.  Composing $\Phi$ with this automorphism, we may assume
WLOG that $\beta = I$, so that $(a,u)\Phi = ( (a)\alpha, u)$.
Then, since $\Phi$ is an isomorphism, we have:
\[
((ab)\alpha,\  (u)\tau_b + v)
= ( (a,u) \cdot_\tau (b,v) ) \Phi = (a,u)\Phi \cdot_\sigma (b,v)\Phi =
((ab)\alpha ,\  (u)\sigma_{(b)\alpha} + v)  \ \ ,
\]
so $\tau = \alpha \sigma$.
\end{proof}

It follows now that the number of nonassociative extra loops of order $16p$ is
independent of $p$.  In the case $p = 3$, that number is already 
known to be $16$, since Goodaire, May, and Raman \cite{GMR}, 
following the classification of Chein \cite{CH2}, 
have listed all nonassociative Moufang loops of order less than $64$.
From Appendix E of \cite{GMR}, we find that $16$ of the Moufang
loops of order $48$ are extra loops.

\begin{theorem}
\thmlabel{sixteenp}
For each odd prime $p$, there are exactly $16$ nonassociative
extra loops of order $16p$.
\end{theorem}


\section{Conclusion}
\seclabel{conclusion}

Although this paper has focused on extra loops, many of the lemmas
hold more generally for CC-loops.
For example, if $Q$ is a CC-loop, then by 
{\cyr Basarab} \cite{BAS}, $Q/N$ is an abelian group.
Of course, $Q/N$ need not be boolean, but if $Q$ is power-associative,
then $Q/N$ has exponent $12$.
Also, if $Q$ is power-associative, nonassociative, and finite,
then $|Q|$ is divisible by either $16$ or $27$.
These results on power-associative CC-loops will appear elsewhere \cite{KK}.


\begin{acknowledgement}
We would like to thank M.~Aschbacher for
suggesting the proof of Lemma \lemref{SylowB},
which is somewhat shorter than our original proof.
\end{acknowledgement}


\end{document}